\documentstyle[12pt,seceqn]{article}
\font\BBB=msym10 scaled \magstep 1
\newtheorem{Theorem}{Theorem}

\title {Ergodic theory for smooth one dimensional dynamical systems}

\author{Misha Lyubich}
\date{October 14, 1990,\hspace{.1in} revised July 2, 1991}

\newfont{\mathsym}{msym10 scaled\magstep1}
\def\ssm{\mbox{\mathsym\symbol{'162}}}

\begin{document}
\maketitle
\begin{abstract}
    In this paper we study measurable dynamics for the widest
reasonable class of smooth one dimensional maps. 
Three principle  decompositions are described
in this class :  decomposition of the global measure-theoretical
attractor into primitive ones, ergodic decomposition and Hopf decomposition.
 For maps with negative Schwarzian 
derivative this was done in the series of papers [BL1-BL5], but the approach
to the general smooth case must be different.
\end{abstract}

\hspace{2in}{\bf Notations}

\vspace{.1in}
$X^\circ\equiv{\rm int} X$ is the interior of a set X;

$\overline X\equiv{\rm cl}X$ is the closure of X;

$[x,y]$ is a (closed) interval ending at $x$ and $y$ (without assuming
         $x\leq y$);

$[U,V]$ is the closed convex hull of sets $U$ and $V$;

$f^n$ is the n-fold iterate of a map $f$;

${\rm orb}(x) = \{f^n x\}_{n=0}^\infty$ is the {\em orbit (trajectory)} of $x$;

${\rm orb_p(x)}=\{f^n x\}_{n=0}^p$;

$\omega (x)$ is the limit set of ${\rm orb}(x)$;

{\BBB N}$=\{1,2,...\}$ is the set of natural numbers.

\section {Statement of the results}
 
In 1985  John Milnor [M] suggested a new approach to measurable dynamics
based upon a concept of (measure-theoretical) attractor. He showed that
any smooth dynamical system has a unique global attractor $A(f)$ and stated the 
problem of  decomposing it into minimal ones. Then the minimal attractors
would give a view of the structure of  typical $\omega$-limit sets.

In the series of papers [BL1-5] the following realization of this program
for one dimensional maps with negative Schwarzian derivative was given
(alternative approaches in the $S$-unimodal case were found in [GJ] and [K]).
The decomposition of $A(f)$ into ``primitive" 
ones was described (slight modification of 
the ``minimality" property was necessary),
and it turned out that it is intimately related to two main measure-theoretical
decompositions: ergodic decomposition and Hopf decomposition.

The aim of this paper is to develop this theory under  proper
smoothness assumptions (without the negative Schwarzian derivative
condition).

Let $M$ be a closed interval,
 and  $\cal A$ denote a class of $C^2-$smooth maps  $f: M\rightarrow M$
with finitely many non-flat critical points (see \S 2 for the accurate 
definition).
 Denote by $\lambda$ the Lebesgue measure on $M$. 
For an invariant set $A\subset M$ let 
\[{\rm rl}(A)=\{x\in M : \omega(x)\subset A \}\]
\[{\rm RL}(A)=\{x\in M : \omega(x)= A \}.\]
These are two ways of understanding of {\em the realm of attraction of $A$};
 we need both of them.

By a {\em limit cycle}  we mean a periodic orbit $Z=\{f^k a\}_{k=o}^p$
whose realm of attraction has non-vacuous interior. 
An interval $I$ is called {\em periodic of period p} if $f^p I\subset I$.
Then
\[{\rm orb}_p (I)=\bigcup_{m=0}^p f^m (I)\]
is called a {\em cycle of} $I$.
If besides $f^p|I$ is monotone, $I$ is called {\em a periodic homterval}.
Any orbit originating in such an interval converges to a cycle,
but it can happen that the set of non-limit cycles in orb$I$ has positive 
measure. This circumstance forces us to take care of periodic homtervals. 
Set
\[\Lambda(f) = M\ssm(\cup\;{\rm  rl}^\circ (Z_i)
\bigcup(\cup_{n=1}^\infty f^{-n}O_i^\circ)),\]
the union is taken over all limit cycles $Z_i$  and all cycles $O_i$ of 
periodic homtervals. What we have removed from $M$ is the trivial part of the
dynamics.

\begin{Theorem}
  The restriction $f\mid\Lambda(f)$ has only finitely many  ergodic components 
  $E_i$.
\end{Theorem}

This means that there is the decomposition 
$\Lambda_f =\cup_{i=1}^t E_i \,(\rm {mod\, 0})$
of $\Lambda_f$ into the finite union of invariant sets of positive
measure such that $f\mid E_i$ are ergodic (where ``mod 0" means that we
ignore sets of measure zero).

Now let us introduce an important notion of {\em attractor} in the sense
of Milnor [M]. This means a closed invariant set $A\subset M$ such that

\noindent
(i) $\lambda(\rm rl (A))>0$

\noindent
(ii) $\lambda(\rm rl(A)\verb+\+\rm rl (A'))>0$ for any proper closed invariant
subset $A'\subset A$.

It is shown in [M] that there is a unique {\em global attractor} $A(f)$, i.e.
an attractor for which ${\rm rl} (A(f))=M$ (mod 0).
 The same is true for the restriction
of $f$ onto any closed subset $K\subset M$ . 
The corresponding global attractor
will be denoted by $A(f|K)$ (if $\lambda(K)=0$ then $A=\emptyset$ ).

Clearly, if $\lambda({\rm RL}(A))>0$ then $A$ is an attractor. 
Such attractors we call {\em primitive}. The ergodic
decomposition will allow us to obtain a decomposition of the global attractor
into primitive ones. Let us say that an orb($x$) is {\em absorbed} by an invariant
set $O$ if $f^n x\in O$ for some $n\in${\BBB N}.

\begin{Theorem}
There is a decomposition 
\begin{equation}
A(f|\Lambda(f))=\cup A_k
\end{equation} 
of the global attractor into the union of
finitely  many infinite primitive attractors $A_k\subset\Lambda(f) $. 
Moreover,

\vspace{2pt}\noindent
({\rm i}) For almost all $x\in M$ either orb$(x)$ is absorbed by 
   a cycle $O_i$ of periodic homtervals,
          or it tends to a limit cycle $Z_j$, or $\omega(x)=A_k$ for some $k$;

\vspace{3pt}\noindent
({\rm ii}) Each $A_k$ contains a critical point; 

\vspace{3pt}\noindent
({\rm iii}) The intersection of any two primitive attractors
            is at most finite; 

\vspace{3pt}\noindent
({\rm iv}) ${\rm RL}(A_k)={\rm rl}(A_k)\cap\Lambda (f)=E_k$ (mod 0) 
for some ergodic component $E_k$. This gives one-to-one correspondence
between primitive attractors $A_k$ and ergodic components 
$E_k\subset \Lambda(f)$
\end{Theorem}
{\bf Remark 1.1. } It is proved in [MMS] that limit cycles of $f\in {\cal A}$
have uniformly bounded periods. Consequently, for {\em analytic} $f$ 
the whole number of limit cycles $Z_k $ and cycles $O_j$ of periodic homtervals
is finite (note that each $O_j$ either contains a limit cycle or
$f^p|O_j\equiv$id).
Moreover, if $M$ is an interval
and $f^2\neq$ id (or $M$ is a circle and $f^p\neq$ id for any $p$ ), 
then  $A(f)$ can be decomposed
into the union of finitely many primitive attractors:
\[A(f)=\cup Z_j\cup A_k.\]

A transformation $f$ of a measure space $(X,\mu)$ is called {\em conservative}
if it satisfies the conclusion of the Poincar\'{e} Return Theorem:
for any measurable subset $Y\subset X$ almost all points $y\in Y$ return to
$Y$ infinitely many times.
The {\em conservative kernel}  $K(f)$
 of $f$ is a maximal measurable set such that
$f\mid K(f)$ is conservative. By {\em Hopf decomposition} one means the
decomposition of $M$ into conservative part $K(f)$ and dissipative part
$M\ssm K(f)$. The following result shows that in one dimensional situation 
it can be described through the notion of attractor. 

\begin{Theorem}
The global attractor $A(f)$ coincides mod 0 with the conservative
  kernel $K(f)$.
\end{Theorem}

An invariant set $K\subset M$ is said to be {\em topologically minimal}
if $\omega(x)=K$ for {\em all} $x\in K$. By a {\em Cantor attractor}
we mean an attractor which is a Cantor set.

\begin{Theorem}
Let $A$ be a Cantor primitive attractor. Then

\noindent
({\rm i}) the restriction $f\mid A$ is topologically minimal;

\noindent
({\rm ii}) topological entropy $h(f\mid A)$ is equal to zero;

\noindent
({\rm iii}) there is a critical point $c\in A$ such that $A=\omega(c)$.

\noindent
({\rm iv}) $A$ does not intersect any other sets of the decomposition (1.1).
\end{Theorem}
{\bf Remark 1.2.} Till now it is unknown if there are Cantor attractors
  different from Feigenbaum-like ones (see the next section for the definition).

\begin{Theorem}
Let $A$ be a primitive attractor. Then one of the following possibilities holds:

\noindent
{\rm A1}. $A$ is a limit cycle;

\noindent
{\rm A2}. $A$ is a cycle of transitive intervals;

\noindent
{\rm A3}. $A$ is a Cantor attractor.
\end{Theorem}

\vspace{.1in}\noindent
{\bf Corollary 1.1.}{\em The number of infinite primitive attractors $A_k$
(which is equal to the number of ergodic components $E_k\subset\Lambda(f)$)
does not exceed the number of critical points in $\Lambda(f)$.}

\vspace{.1in}
\noindent
Let us say that $A$ is a {\em minimal} attractor  if there are no smaller
attractors $A'\subset A$.  

\vspace{.1in}\noindent
{\bf Corollary 1.2.}
{\em Any primitive attractor $A$ is minimal except only one case:  $A$ is
a cycle of transitive intervals whose boundary $\partial A$  contains a
 parabolic limit cycle. }

\vspace{.1in}
\noindent{\bf Remark 1.3.}
The whole above theory holds for maps of the circle as well but the 
statements need minor modifications concerning  
 immersions of the circle. Let us mention also that 
ergodicity of circle diffeomorphisms was proved quite long ago by
M. Herman [H] and A. Katok (see [KSF]).

Let us describe the structure of the paper. Section 2 contains preliminaries
on topological one dimensional dynamics (including the principle concept
of a chain of intervals) and distortion lemmas.
Sections 3,4 are technical ones. The former contains an estimate of the
intersection multiplicity of a monotone chain of intervals.
The latter
explains how to control density moving along a chain of intervals
(it needs a concept of ``$D(X,\epsilon)$ broken lines").
 Getting together all these tools,
in the last Section 5 we prove the main results.
    
 This paper is a natural continuation of the series [BL1-BL5].
However, the approach here differs from that in previous ones . 
In the case of negative Schwarzian derivative we started from the description
of primitive attractors , and then pass to ergodic decomposition etc. 
In the general smooth case 
considered here the way is opposite (and it is the only way known to the 
author): the starting point is ergodic decomposition while the primitive
attractors can be described only in the very end.

Finally, I would like to thank  A. Blokh and J. Milnor 
for looking through the manuscript and making  useful comments.

\section{Preliminaries}
 
Let $M$ be a finite union of disjoint closed intervals (it is convenient
for technical reasons to consider non-connected $M$ as well).
Remember that ${\cal A}$ denotes a class of $C^2$-smooth self-maps
of the manifold $M$ with non-flat critical points. The latter means that
there are $C^2-charts$ around each critical point  $c$ and the 
corresponding critical value
$fc$ in which $f$ is reduced to the form
$x\mapsto\sigma |x|^r$ with real $r\geq 2$ and
a sign $\sigma\in \{+1,-1\}$
which may depend on the sign($x-c$).  

\vspace{.1in}
{\bf Remark 2.1.} The smoothness conditions on the map $f$ are determined
 by the range of 
validity of so called Koebe Principle (see below). The precise 
regularity for it ($C^{1+Zygmund}$) was established by
Sullivan [S]. So, this seems to be the widest reasonable regularity for our
theory.
\vspace{.1in}
  
Denote by $C(f)$ the set of critical points of $f$. There are critical points
of two types: {\em extrema} or {(\em turning points)}  and
 {\em inflection points}. Let $d$ denote the number of extrema.
Points of the set $C(f)\cup\partial M$ will be called {\em singular}.

There is a natural involution $\tau$ in a neighbourhood of any extremum
 $c$, namely
$\tau (x)=x'$ if $f(x)=f(x')$. It follows from non-flatness that $\tau$ is
smooth with $\tau '(c)=-1$.

Actually, by $C^2$-smooth conjugacy the map $f\in {\cal A}$ can be reduced to
$\sigma |x-c|^r+b$ in a neighbourhood of any critical point $c$.
{\em In what follows we will suppose that it is the case.}
Then $\tau$ is reduced
to the standard isometric reflection with respect to $c$.

Let us make also the following convention:
\begin{equation}
f(\partial M)\subset\partial M .
\end{equation}
It is possible because of the following {\em surgery}. Let us include $M$ into
a compact one dimensional manifold $\tilde M$ such that 
$\partial\tilde M\cap M =\emptyset$. Then $f$ can be continued to a map
$\tilde f\in {\cal A}$ of  $\tilde M$
in such a way that $C(\tilde f)=C(f)$ and (2.1) holds for 
$\tilde f :\tilde M\rightarrow\tilde M$.

 Now we need more definitions.
An invariant closed set $R\subset M$will be called {\em transitive} if it contains a dense orbit.

A set $R$ will be called {\em a basic set} (see [B1, B2])
 if it is a set of all points belonging to a cycle 
of intervals $ O=orb(I)$ of period $p$
 and satisfying the following property:
for any open interval $J\subset O$
intersecting $R$ and for any compact subset $K\subset {\rm int}O$
there is a $N$ such that 
\begin{equation}
\bigcup_{k=n}^{n+p} f^k J\supset K
\end{equation} 
 for all $n\geq N$. Clearly, a basic set is closed and invariant.
There are basic sets of three types : periodic orbits, cycles of intervals 
and  Cantor basic sets.

By {\em a Feigenbaum-like attractor} we mean an invariant Cantor set $F\subset M$
of the following structure:
\[ F=\bigcap_{n=1}^\infty orb_{p_n} (I_n)\]
where $I_1\supset I_2\supset ...$ is a nested sequence of periodic intervals
of periods $p_n$ such that $p_n\rightarrow\infty$.

The following Topological Structural Theorem
 follows from the pure topological considerations
(see [Sh, JR , Ho , B1 , B2 ]) and the absence of wandering intervals
(see  [G1], [Y], [L], [BL7], [MMS] and references there). 

\vspace{.2in}
\noindent
{\bf Theorem A.} 
{\em For any $x\in M$ one of the following possibilities holds:

\noindent
{\rm (0)} ${\rm orb}(x)$ is absorbed by a cycle of a periodic
           homterval;

\noindent
{\rm (i)} ${\rm orb}(x)$ tends to a limit cycle;

\noindent
{\rm (ii)} ${\rm orb}(x)$ is absorbed by a basic set;

\noindent
{\rm (iii)} $\omega (x)$ is a Feigenbaum-like attractor.}

\vspace{.2in}
Now let $R$ be either a basic set or a Feigenbaum-like attractor. 
Restricting  $f$ onto an appropriate cycle of intervals
and using the above surgery we can {\em localize $f$ with respect to $R$} 
in the following sense:

\noindent
{\rm L0.} $R\subset{\rm int}M.$

\noindent
{\rm L1.} All critical points belong to $R$.

\noindent
{\rm L2.} There are no limit cycles in int$M$.

Maps satisfying L1-L2 we will call $R$-{\em local} (or just {\em local}).
 There are local maps of
two types: {\em finitely renormalizable} when $R$ is a basic set,
and {\em infinitely renormalizable} when $R$ is a Feigenbaum-like attractor.

It is easy to see that the mixing property (2.2) together with L2 yield
{\em the sensitive dependence to initial conditions} on a local basic
set $R$ in the following sense.  There exists $\gamma(R)>0$ with the following
property: $\forall\tau \,\exists N$ such that for any {\em closed} interval $J$
intersecting $R$, $\lambda(J)\geq\tau$,
 we have $\lambda(f^n J)>\gamma(R),\; n\geq N.$

The following easy but useful Proposition was stated in [L].

\vspace{.1in}\noindent
{\bf Proposition 2.1.}
{\em Let $J$ be an interval whose orbit does not tend to a limit cycle. Then}
\[\inf_{0<m<\infty}\lambda(f^m (J))>0.\]

\vspace{.1in}
Let us fix two constants $\xi >0$ and $\eta >0$ till the end of the paper.
Let $\eta$ be so small that $\eta$-neighbourhoods of critical points don't
intersect, and the involution $\tau$ is well-defined in the $\eta$-neighborhoods
of extrema.
Then choose an $\xi$ by Proposition 2.1 in such a way that
for any interval $J$ containing a critical point
\begin {equation}
 \lambda (J)>\eta \Rightarrow \lambda (f^m (J)) >\xi 
\end{equation}

Now let us introduce  a notion of {\em a maximal chain of intervals} which is a 
key to polymodal maps (see [L]). By {\em a chain of intervals}  {\BBB I}
we mean just a sequence of intervals $\{I_m\}_{m=1}^n$ such that
$f I_m\subset I_{m+1}, m=0,1,...,n$. The chain is called {\em maximal}
if $I_m$ are the maximal intervals satisfying this property.

Let $I_n\subset M^\circ$ and $\lambda (I_n)<\xi$. Then we have for the maximal
chain {\BBB I}  and $1\leq m < n$ that (take into account (2.1))
either $f: I_m\rightarrow I_{m+1}$ is a homeomorphism,
or $I_m$ is symmetric with respect to some extremum, and 
$f(\partial I_m)\subset\partial I_{m+1}$

In what follows $I$ denotes a closed interval such that $I\subset M^\circ$
and $\lambda(I)<\xi$.

Let  $x\in M$ and $f^n (x)\in I$.
The main way of constructing  maximal chains of intervals is a 
{\em pull-back} of  $I$ along the ${\rm orb }_n(x)$ .
Namely, set $I_n=I$, and $I_m$ be the maximal interval containing $x_m$
for which $f^{n-m}I_m\subset I,\, m=0,1,...,n-1.$

Define the {\em order} ord\,{\BBB I}$\equiv$ord$(n,x,I)$
 of the chain (pull-back) {\BBB I} as the number of intervals
$I_m$ containing extrema.
If $f^n$ monotonously map $I_0$ onto $I_n$ we say that {\BBB I} is a 
{\em monotone} chain (pull-back).

\vspace{.1in}
{\bf Proposition 2.2} [L]. Consider a local map $f$.
Let $n$ be the first moment when orb($x$) passes through $I$ .
Consider a pull-back {\BBB I}=$\{I_m\}_{m=0}^n$ of $I$ along orb$_n(x)$, and
 let $\{I_{m_i}\}_{i=1}^{\nu}$ be the intervals of
the chain  containing extrema. Then we have for $i>d$
(where $d$ is the number of extrema)\newline
(i) $I_{m_{i-d}}\supset I_{m_{i}}$ and $I_{m_i}$
 is periodic;\newline
(ii) $m_i$ is the first moment when orb$(x)$ passes through int$I_{m_{i-d}}$.
\newline
(iii) $f^{m_i}$ monotonously maps a neighbourhood of $x$ onto an appropriate
half of the interval $I_{m_{i-d}}$.

\vspace{.1in}
{\bf Corollary 2.1.} In the above situation provided $I$ is non-periodic,
we have ord\,{\BBB I}$\leq d$. In particular, it is the case if $f$ is
finitely renormalizable and $I$ is small enough.

\vspace{.1in}
Now let us describe the analytical tools of the paper in the form 
of two Distortion Lemmas. Both of them follow from {\em the Koebe Principle}
in one dimensional dynamics intensively exploited in recent works, see
[Y], [G2], [MS], [L], [Sw], [S] ...

In what follows we assume that $f\in{\cal A}$, $f^n$ is {\em monotone}
(perhaps, with critical points)  on an interval $J$
and denote by {\BBB J} the monotone chain
of intervals $\{f^l J\}_{l=0}^n$.

Denote by $\mu=$mult{\BBB J} the {\em intersection multiplicity} of {\BBB J}
(i.e. the maximal number of intervals from {\BBB J} with non-empty
intersection).
For a measurable set $X$ let

\[{\rm dens}(X|I)=\lambda(X\cap I)/\lambda(I),\]
\[{\rm Dens}_a(X|I)\equiv{\rm Dens}(x|[a,b])=
  \sup_{y\in(a,b]} {\rm dens}(X|[a,y]).\]
In what follows assume that $X\subset I$.

\vspace{.1in}
{\bf Three Interval Distortion Lemma} (see [MS], [BL6]).
{\em  Consider an interval $I\subset J^\circ$,
and let $J^+$ and $J^-$ be the components of $J\ssm I$. 
Suppose \[\lambda(f^n J^{\stackrel{+}{-}}/\lambda(f^n I)
\geq\delta.\]
 Then there are  positive constants
 $\sigma=\sigma_\mu(\delta)$  and 
$q(\epsilon)=q_{\mu}(\epsilon, \delta)$
, $q(\epsilon)\rightarrow 0$ as $\epsilon\rightarrow 0$, such that}
\newline (i) $\lambda(J^{\stackrel{+}{-}})/\lambda(I)
\geq\sigma_\mu(\delta).$
\newline (ii) dens$(f^n X|f^n I)\geq 1-\epsilon\Rightarrow$
   dens$(X|I)\geq 1-q(\epsilon)$.

\vspace{.1in}
{\bf Two Interval Distortion Lemma.} (see [BL6]).
{\em  Divide $J$ by a point $a$ into two subintervals $L$ and $R$,
$ b=f^n a$. Assume 
$\lambda(f^n L)/\lambda (f^n R)\leq~K$. Then
\[{\rm Dens}_a(X|L)\geq 1-\delta\Rightarrow
{\rm Dens}_b(f^n X|f^n L)\geq 1-\alpha_\mu(\delta, K)\]
where $\alpha_\mu(\delta, K)\rightarrow 0$ as $\delta\rightarrow 0$, $K$ fixed.}

\section{An estimate of the intersection\newline multiplicity}

A technique of estimating the intersection multiplicity of a monotone
pull-back of an interval $I$ was developed by A. Blokh (see [BL7], \S 2.3)
in order to generalize the results of [L] onto the smooth case.
The interval $I$ was supposed to be symmetric around an extremum. 
Here we will develop the technique for an arbitrary $I$ (concentrating
only over new points).  

Let us pass to the  main definition.
Consider an interval $I=[a,b]$ and a point $x\in M$ which
 does not lie in a basin of a limit cycle.
  Let $x(n)=f^n x$ be the first point of orb$(x)$ 
lying in $(a,b)$. Assume that there are $p,r\in${\BBB N} and a point $v\in (a,b)$
for  $m=n-(r-1)p$ and 
$J=[x(m), v]$  the following properties hold\newline 
D1. $J^\circ$ contains $x(n)$ and exactly one of the points $a, b$, say $a$;
\newline 
D2. $f^p$ is an orientation preserving homeomorphism of $J$ onto
   $[x(m+p), b]$;\newline
D3. $J^\circ$ contains points $x(n+ip),\; |i|\leq r-2$, and no other points
   of the orb$_{n+(r-2)p}x$.   

Then we say that $x(n)$ belongs to a {\em multiple collection}
 $\{x(n+ip)\}_{i=-(r-1)}^{r-2}$. A number $r$ is called the
{\em depth} of $x(n)$ in the collection. 
Let us denote by {\em depth} dp$_a(n,x,I)\equiv$dp$_a(n)$ 
the maximal depth of $x(n)$
in all multiple collections containing  it (if there are no such collections
set dp$_a(n)=0$). Set

\begin{equation}
{\rm dp}(n)\equiv {\rm dp}(x,n,I)= \max\{{\rm dp}_a(n), 
  {\rm dp}_b(n)\}.
\end{equation}

\vspace{.1in}
{\bf Remark 3.1.} Observe that if dp$_a(n)>1$ then $f^p$ move all points of $J$
 toward $b$. Indeed, it is true for  the endpoint $x(m)$. If it fails
for some point of $J$ then $J$ contains a fixed point $\alpha$. Hence, 
$f^p$ maps monotonously the interval $[x(m),\alpha]$ into itself.
Then orb$(x)$ should converge to a cycle contradicting the assumption.   

\vspace{.1in}
{\bf Remark 3.2.}  Observe also that
$\min\{{\rm dp}_a(n), {\rm dp}_b(n)\}=1$.
Indeed, denote by $r_a, p_a, J_a$ the data corresponding to $a$,
and use the similar notations for $b$. Assume that $r_b\geq r_a>1.$
Then $x(n+p_a)\in [a,x(n)]\subset J_b$ contradicting D3.

\vspace{.1in}
{\bf Lemma 3.1.} {\em Under circumstances described above
let H be an interval ending at x such that $f^n$
monotonously maps H onto $[x(n), b]$. Assume ${\rm dp}(n)\geq 2$. Then}

\[{\rm mult}\{f^k H\}_{k=0}^n \leq 2{\rm dp}_a(n).\]

{\bf Proof.} Let a point $y$ belong to $\kappa$ of the intervals 
$H_\kappa\equiv f^\kappa H$, $r$ of them lying on the one side of $y$
and $\kappa-r$ on the other. We are going to prove that $r\leq {\rm dp}(n)$
and $\kappa-r\leq {\rm dp}(n)$, which implies the required. Clearly, we
can restrict ourselves to the estimate of $r$.

Let $[x(i_1), y]\supset...\supset[x(i_r),y],\; [x(i_k), y]\subset H_{i_k}$ .
Then we have (see [BL7])\newline
(i) $i_1<...<i_r$;\newline
(ii) $x(i_k)$ are the only points of the orb$_n (x)$ lying in 
    $K\equiv [x(i_1), x(i_r)]$.

Applying $f^{n-i_r}$ we can assume that $i_r=n$. Denote $i_1=m$ .
Assume also for  definiteness that $a<b$.

By the assumption, all ponts $x(i_k)$ lie outside $(a,b)$. In fact, they
should lie to the left of $a$. Indeed, otherwise 
$H_m\supset  [x(m), x(n)]\supset [x(n), b]$. But $f^{n-m}$ maps
monotonously $H_m$ onto $[x(n), b]$. So, orb$(x)$ would have converged to
a limit cycle. 

Now set $p=n-i_{r-1}$, and apply $f^p$ to $K=[x(m), x(n)]$. 
Since $K\subset H_m$ and $p<n-m$, $f^p$ is monotone on $K$. Since
$f^p$ maps $[x(i_{r-1})
, x(n)]\subset H_{i_{r-1}}$ into 
$H_n=[x(n), b]$ preserving orientation, we conclude that
$f^p |K$ preserves orientation. Now it follows from the above property (ii)
that (see [BL7])

\[i_{k+1}=i_k+p , \;k=1,..., r-1.\]

Further, let $v$ be the right endpoint of the interval $H_{n-p}$. Denote
$J=[x(m), v]$. Clearly, $J$ satisfies properties D1-D2 of the definition
of a multiple collection. Let us check that it satisfies D3 as well.

First, $x(n-ip)\in K\subset J$ for $i=1,...,r-1.$ Further, 
$x(n+ip)\in f^{ip}H_{n-ip}=[x(n),b]$ for $i=1,...,r-1$. Moreover, for
$i<(r-1)$ we have $x(n+ip)\in f^{ip}H_{n-(i+1)p}=[x(n-p),v]$, so
$x(n+ip)\in J^\circ$.

Let $x(l)\in J^\circ,\,\; l\leq n+(r-2)p,\; l\neq n+ip$ for $|i|\leq r-2$.
Then $l>n.$ Indeed, $J\subset (a,b)\cup K.$
If $l\leq n$ then $x(l)$ does not belong to $(a,b)$ by the assumption
and does not belong to $K$ by the above property (ii).

So, $l=n+k$ where $0<k<p(r-2)$. Hence, $s\equiv m+k<n-p$, and $x_s$ lies
outside $[x(m), b]$ according to what has been proved right now.
Hence, the interval $f^k K=[x(l), x(s)]$ contains one of the points
$x(m), b$. But it cannot contain $x(m)$ because otherwise
\[(x(n),b)\supset {\rm int}(f^{n-m} K) =
 {\rm int}(f^{n-m-k}(f^k K))\ni x(n-k).\]
Consequently, $T\equiv [x(l), x(s)]\ni b$. Moreover, $f^p$ maps $T$
monotonously, orientation preserving and without
fixed points onto $[x(l+p), x(s+p)]$. Since $x(s+p)\in T$, $f^p$
moves all points of $T$ to the left. In particular,  $b$ is moved to the left.

Further, $f^p$ is also a monotonous map on the interval $[a, v)$ which has
 with $T$ a common point $x(l)$. Hence, $f^p$ is monotonous on $[a,b]$. Since 
 $a$ is moved to the right and $b$ to the left, all orbits in
$[a,b]$ converge to limit cycles contradicting the assumptions.

So, the interval $J$ satisfies all properties $D1-D3$ and hence
dp$_a(n)\geq r$.   

\vspace{.1in}
{\bf Corollary 3.1.} {\em Under the above circumstances let 
{\BBB I}=$\{f^k I\}_{k=0}^n$ be the monotone pull-back of $I=I_n$
along orb$_n(x)$}. Then
\vspace{.1in}

\hspace{1.8in}
mult\,{\BBB I}$\leq 2 ({\rm dp}(n)+1)$

\vspace{.1in}
{\bf Proof.} Divide $I_0$ into two intervals $H^+$ and $H^-$ ending
at $x$, and apply the lemma to these intervals taking into account Remark 3.1.

\section{Transportation of broken lines}
Here we are going to prove the series of density lemmas. 
In what follows $f$ is assumed to be $R$-local, $X$ denotes an invariant
set of positive measure, $I=[a,b]\subset {\rm int} M$.

By a {\em broken line} beginning at $x$ and ending at $y$ we mean a sequence of 
points ${\cal L}=\{x_k\}_{k=0}^n$ such that $x_0=x , x_n=y$; the intervals
$[x_k, x_{k+1}]$ are the {\em links} of the broken line 
(a link can be degenerate, i.e. $x_k=x_{k+1}$). We say that
$\cal L$ is a {\em proper} broken line if all links are non-degenerate
and \[[x_{k-1}, x_k]\subset [x_k, x_{k+1}],\; k=1,...,n-1\]. 

We say that ${\cal L}$ satisfies $D(X,\epsilon)-property$
if for any non-degenerate link $[x_k,x_{k+1}]$ we have
${\rm Dens}(X|[x_k, x_{k+1}])\geq 1-\epsilon$ for $k=1,2,...,n-1$. 
Any $D(X,\epsilon)$-broken line can be easily turned into a proper
$D(X,\epsilon)$-broken line with the same beginning and end .

\vspace{.1in}
{\bf Lemma 4.1.} {\em Let a point $x$ don't
converge to a limit cycle, $n$ be the first moment when {\rm orb}$(x)$
passes through $I^\circ=(a,b)$. Assume that there is an interval $H\ni x$
monotonously mapped onto $I$ by $f^n$.
Let $H^+$ and $H^-$ be the closures of the components of $H\ssm\{x\}$.
Then $\forall\epsilon\; \exists\delta$
such that if} 

\vspace{.1in}
\hspace{.8in}
Dens$_x(X|H^+)>1-\delta$ or Dens$_x(X|H^-)>1-\delta$

\vspace{.1in}
\noindent
{\em then there exists $D(X,\epsilon)$ broken line 
beginning at $x_n\equiv f^n x$ and ending
at an endpoint of $I$.}

\vspace{.1in}
{\bf Proof.} Let us make the following conventions:
dp$_b(n,x,I)=1$ (see Remark 3.1) and $f^n H^-=[a, x_n]$.
Denote $r={\rm dp_a}(n,x,I)$. Then by Corollary 3.1 we have
mult$\{f^l H\}_{l=0}^n\leq 2r+2$.

Let $\epsilon>0$. Let us choose a big number  K 
(so that $(1+K^{-1})(1-\epsilon)<1$);
then choose $\sigma=\sigma_6 (1/2K)$
by the Three Interval Distortion Lemma,
and then find a small number
$\delta$ satisfying the following inequalities:

\[\frac{1-\alpha_6(\delta, 2K)}{1+K^{-1}}>1-\epsilon,\qquad
\alpha_6(\delta, \sigma^{-1})<\epsilon\]
where $\alpha$ is taken from Two Interval Distortion Lemma.

Assume first that $r\leq 2$. If $|a-x_n|\leq K |b-x_n|$
then by the Two Interval Distortion Lemma we get
Dens$(X|[x_n,a])\geq 1-\alpha_6 (\delta, K)> 1-\epsilon.$ 

Otherwise consider an interval $[y,x_n]\subset [a,x_n]$ such that 
\[|y-x_n|=K|b-x_n|.\]
Then Dens$(X|[x_n,y])> 1-\epsilon$, and for any $w\in [x_n, b]$ 

\begin{equation}
{\rm dens}(X|[y,w])\geq  {\rm dens}(X|[y, x_n])\frac{|y-x_n|}{|y-b|}
\geq \frac{1-\alpha_6(\delta, K)}{1+K^{-1}}> 1-\epsilon.
\end{equation}
It follows  that the two-linked broken line
$\{x_n,y,b\}$ can be turned
into a  $D(X,\epsilon)$ broken line $\{x_n,\overline y,b\}$ with
$\overline y\in [y,x_n]$ (perhaps, $\overline y=x_n$). To this end
it is enough to set $\overline y$ to be the nearest to $x_n$ point of
$[y,x_n]$ satisfying (4.1).

Assume from now on that
 $r\geq 3$, and set 
\[w_i=x(n+ip)\in I,\; i=0,1,...,r-2,\; w_{-1}=a,\; w_{r}=b.\]
Denote by $M_i$ the intervals on which these points divide $I,\; i=0,1,...,r-1$.
In particular, $ M_0=[a,x(n)], M_{r-1}=[w_{r-1}, b]$. 
Now consider two cases:

\vspace{.1in}
I. Dens$_x(X|H^-)>1-\delta$. Then let us consider two subcases:

\vspace{.1in}
(i) $\lambda(M_{k+1})\geq K^{-1}\lambda(M_k)$ for some $k\in[0,r-2]$.
Let us take the first such  $k$.
Assume $k=0$. Then consider the interval $G\subset H$ mapped
onto $M_0\cup M_1$ by $f^n$.  It is easy to see that
dp$(n,x,M_0\cup M_1)=2$.
Applying  the Two Interval Distortion Lemma to $f^n|G$
we get

\[{\rm Dens}(X|[x(n), a])> 1-\alpha_6 (\delta,K)>1-\epsilon.\] 
So, $\{x(n), a\}$ is a $D(X,\epsilon)$ one-linked broken line.

Now let $k>1$.  Then we have
\[\lambda(M_{k+1})\geq\frac{1}{K}\lambda(M_k) \; {\rm and}\;
 \lambda(M_{k-1})\geq\lambda(M_k). \]
Let $k<r-2$.
Applying the Three Interval Distortion Lemma to $f^{kp}$
with the central interval $M_0$, 
 we get $\lambda(M_1)\geq\sigma\lambda(M_0)$. Now we can apply the Two Interval
Distortion Lemma to $f^n$ as above. It gives

\[{\rm Dens}(X|\,[x(n), a])> 1-\alpha_6(\delta,\sigma^{-1})>1-\epsilon,\]
and we are done.

Finally, let $k=r-2$ . Then consider an interval 
$N_1=[w_0,z],\; w_1\leq z< w_2,$
monotonously mapped onto $M_{r-2}$ by $f^{(r-2)p}$. Replacing $M_1$ by $N_1$
in the above argument we will get  the same conclusion. 

\vspace{.1in}
(ii) Now assume $\lambda(M_{k+1})\leq K^{-1}\lambda(M_{k})$
 for all $k\in [0,r-2].$
Then for any such $k$ we can construct a two-linked  $D(X,\epsilon)$ broken line
beginning  at $w_{k}$ and ending at $w_{k+1}$ (see the above argument for
$r\leq 2$). Getting together these lines
we obtain a $D(X,\epsilon)$ broken line beginning at $x(n)$ and ending 
at~$b$.

\vspace{.1in}
II. Dens$_x(X|H^+)>1-\delta$. Again let us consider two subcases:

\vspace{.1in}
(i) $\lambda(M_{k+1})\geq 2K\lambda(M_k)$ for some $k\in [0,r-2]$.

If $k=0$ then we have a two-linked $D(X,\epsilon)$ broken line 
in $M_0\cup M_1$ beginning
at $x(n)$ and ending at $a$ (see the argument for $r\leq 2$).

Otherwise consider the first moment $k$ for which 
$\lambda(M_{k+1})\geq 2K\lambda(M_k)$ .
By the same reason as above we have a two-linked
 $D(f^{n+kp}(X\cap H),\epsilon)$
broken line in $M_k\cup M_{k+1}$ beginning at $w_{k}$ and ending at
$w_{k-1}$, and such that
 the length of the first link of the line
does not exceed $K\lambda(M_{k+1})$. So, this line lies deeply
inside the interval $[M_{i-1}, M_{i+1}]$. Hence we can pull it  back 
to $w_0$ by $f^{kp}$   with bounded
distortion (it needs more careful selection of the constants which we leave 
to the reader), and then act as in case I.

\vspace{.1in}
(ii) Assume $\lambda(M_{k+1})\leq 2K\lambda(M_{k})$ for all $k\in [0,r-2].$  

\noindent
Then Dens$(X|\,[w_k, w_{k+1}])>1-\epsilon$  for all $k=0,...,r-2$, 
and hence Dens$(X|\,[x(n), b])>1-\epsilon$. The Lemma is proved.$~~\Box$

\vspace{.1in}
{\bf Lemma 4.2.} {\em Let a point $x\in M$ don't converge to a cycle, 
$n$ be the first moment for which $f^n x\in I^\circ$.
Consider the pull-back {\BBB I}=$\{I_m\}_{m=o}^n$ of $I=I_n$ along ${\rm orb}(x)$,
$\nu ={\rm ord}$\,{\BBB I}. 
\newline
Then $\forall\epsilon >0\;
\exists \delta =\delta (\nu , \epsilon)>0$ such that 
if there is a $D(X,\delta)$ broken line beginning at $x$
and ending at $\partial I_0$ then there is a $D(X,\epsilon)$ broken line
beginning at $f^n x$ and ending at $\partial I$.}

\vspace{.1in}
{\bf Remark 4.1.} If $I$ is non-periodic then by Proposition 2.2 we can 
use this lemma
with $\nu\leq d$. In particular it is the case when $x$ belongs to a
basic set $R$ and $I$ is small enough.

\vspace{.1in}
{\bf Proof.} {\em Step 1}.
  First assume that $\nu$=ord\,{\BBB I}=0, so
$f^n$ monotonously maps $H\equiv I_0$ onto $I$. Let ${\cal L}=
\{x_0=x,...,x_{m-1},x_m\}$ be a given $D(X,\delta)$ broken line
beginning at $x$ and ending at an endpoint $x_m\in\partial H$ of $H$.
Let us check the required by induction in $m$. 
The base of induction $m=1$ is given by Lemma 4.1. 
Without loss of generality we can consider that ${\cal L}$ is  proper.
Then the $m-1$ - linked broken line ${\cal T}=\{x_0,...,x_{m-1}\}$ is contained
in the interval $[a, x_{m-1}]$ where $a\in \partial H$. By the induction
assumption, there exists a $D(X,\epsilon)$ broken line $\cal R$
beginning at $f^n x$
and ending at either $\partial I$ or $f^n x_{m-1}$. In the former case 
we are done. In the latter case construct by Lemma 4.1 a $D(X,\epsilon)$
broken line ${\cal R}'$ beginning at $f^n x_{m-1}$ and ending at $\partial I$.
Getting together ${\cal R}$ and ${\cal R}'$ we obtain the required broken line.

\vspace{.1in}
{\em Step 2}. Let $c$ be an extremum, $J$ be a short $c$-symmetric interval
 containing a 
$D(X,\epsilon)$ broken line ${\cal T}$ . Then ${\cal T}$ can be reconstructed 
into a $D(X, \rho(\epsilon))$ broken line ${\cal T}' $
with the same beginning and end and containing in cl($J\verb+\+ \{c\}$).
Moreover, $\rho (\epsilon)\rightarrow 0 \; {\rm as}\;\epsilon\rightarrow 0$.
The reconstruction is described (implicitly) in [BL4].

\vspace{.1in}
{\em Sten 3}. Consider the intervals $I_{n(1)},...,I_{n(\nu)} $
of the chain {\BBB I} containing extrema. By Step 1 we have a
$D(X,\epsilon_1)$ broken line in $I_{n(1)}$ beginning at $f^{n(1)}$ and
ending at $\partial I_{n(1)}$. By Step 2 we can change it to
$D(X,\epsilon_1')$ broken line whose interior does not contain $c$.
Now apply $f$ to this line using obvious local estimates in a neighborhood
of $c$ . Then we will get a $D(X,\epsilon_1'')$ broken line
in $I_{n(1)+1}$. 

Proceeding in the same manner from $n(1)$ to $n(2)$ , from $n(2)$ to $n(3)$
etc., we will get the required broken line.   $~~\Box$

\vspace{.1in}
{\bf Lemma 4.3.} {\em Let $f$ be finitely renormalizable.
 Let $x$ be a density point of an invariant set $X$
 absorbed by a basic set $R$ , $\epsilon >0$. Then there exist
$\gamma>0$ and a natural $N$
with the following property. For any 
interval $I$ of length $<\gamma$ omitting $N$ subsequent
points of {\em orb}(x) there is an interval $J\subset I\ssm\omega(x)$
such that {\em dens}$(X|L_i)\geq 1-\epsilon$ for each component $L_i$ of
$I\ssm J ,\; i=1,2 $.}

\vspace{.1in}
{\bf Proof. } Clearly, we can assume that  and $x\in R$.  
Now let us select several constants using notations fixed in \S 2.
Let $\lambda (I)<\min\{\eta , \gamma(R) \}\equiv\gamma$.

Choose $\delta_1=\delta(d,\epsilon)$ and $\delta_2=\delta(d,\delta_1)$
by Lemma 4.2 (recall that $d$ is the number of extrema).
Since $x$ is a density point of $X$, there is $\rho>0$ such that
\[{\rm Dens}(X|[x,a])>1-\delta_2\]
if $|x-a|\leq\rho$.
By Proposition 2.1, there is a 
$\tau>0$ such that for any interval $T\ni x$ of
length $>\rho$ we have  
\[\lambda(f^n(T))>\tau,\; n=0,1,2...\]

Finally, by sensitive dependence to initial conditions , 
there exists an $N$ such that for any closed interval
$J$ intersecting $R$ with $\lambda (J) \geq\tau$ we have: 
${\rm diam}\, f^n(J)>\gamma$ for all $n\geq N$. 

Denote $x_m =f^m x$.
Let $I$ omit $N$  subsequent points of ${\rm orb}(x)$ beginning with $x_n$,
and $l\geq N$ be the first moment for which $x_{n+l}\in I$.
Consider the pull-back {\BBB I}=$\{I_m\}_{m=n}^{n+l}$ of $I$ along 
${\rm orb}_{n+l}(x_n)$, \,$I_n\equiv K$.
By the choice of $N$, we have 

\begin{equation}
\lambda (K)<\tau.
\end{equation}

Our nearest goal is to construct a $D(X,\delta_1)$-broken line beginning at
$x_n$ and ending at $\partial K$.
Set \[{\cal P}=({\rm orb}_{n-1}(x)\cap K^\circ)\cup\partial K.\] 
Let $T\subset K$ be the smallest interval containing $x_n$ and ending at
points of ${\cal P}(n)$.

 Consider the  the pull-back $\{T_m\}_{m=0}^n$ of $T$ along ${\rm orb}_n(x)$.
It follows from (4.2) that $\lambda (T_0)<\rho$.
Hence,  $X$ is thick in $T_0\equiv [a,b]$:

\[{\rm Dens}(X|\,[x,a])>1-\delta_2, \;{\rm Dens}(X|\,[x,b])>1-\delta_2.\]
So, the intervals $[x,a]$ and $[x,b]$ can be considered as $D(X,\delta)$-broken
lines (with one links). By Lemma 4.2 , there is a $D(X,\delta_1)$-broken line 
${\cal L}_0$ beginning at $x_n$ and ending at $\partial T$.

If ${\cal L}_0$ is ended at $\partial K$, we are done. Otherwise it ends at
a point $x_{n(1)}\in K$ with $n(1)<n$. Handling  $x_{n(1)}$ 
in the same manner, we
will find a $D(X,\delta_1)$-broken line beginning at $x_{n(1)}$ and ending 
at either $\partial K$ or $x_{n(2)}$ with $n(2)<n(1)$.

Proceeding in such a manner, we construct a sequence of $D(X,\delta_1)$-broken
lines ${\cal L}_0,..., {\cal L}_k$ such that 
 ${\cal L}_{i+1}$ starts at an endpoint of
${\cal L}_i,\;i=0,1,...,k-1,\; {\cal L}_0$ starts at $x(n)$
and ${\cal L}_k$ ends at $\partial K$. Putting together these lines,
we get the required broken line ${\cal L}=\cup{\cal L}_i$.

Now let us consider the map $f^l : K\rightarrow I$. By Lemma 4.2,
we get a $D(X,\epsilon)$-broken line ${\cal Y}$ beginning at $x_{n+l}$ and
ending at $\partial I$. Reconstruct it into a proper $D(X,\epsilon)$-broken
line and consider its last link $S$. This interval contains $x_{n+l}$,
ends at $\partial I$ and ${\rm dens}(X|S)>1-\epsilon$. 

Replace the interval $I$ by $I_1=I\ssm S$ and consider the first moment
$l_1>n$ when $x_{n+l(1)}\in I_1$. Clearly, $l_1>l$. Repeating the previous
argument, we get an interval $S_1\subset I_1$ containing $x_{n+l(1)}$,
ending at $\partial I_1$ and such that dens$(X|I_1)>1-\epsilon$. 

Set $I_2=I_1\verb+\+S_1$ and proceed in the same manner.
 We result with a nested
sequence of intervals $I\supset I_1\supset I_2\supset...$ whose intersection
$J$ satisfies the required properties.     $~~\Box$  

\vspace{.1in}
{\bf Corollary 4.1.} {\em Let $f$ be finitely renormalizable.
Let $I$ be an interval centered at a point $a\in\omega(x)$,
$I^+$ and $I^-$ be the components of $I\ssm \{a\}$.Then }
\[\max({\rm dens}(X|I^+), {\rm dens}(X|I^-))\rightarrow 1 \;{\rm as}\;
\lambda(I)\rightarrow 0.\] 

\vspace{.1in}
{\bf Lemma 4.4.} {\em Let $I$ be a non-periodic interval symmetric around 
an extremum $c$, the set $X$ be $\tau$-symmetric. Then}
$\forall\epsilon\;\exists\delta$ such that

\[\lambda(I)<\delta \Rightarrow{\rm dens}(X|I)>1-\epsilon.\]

\vspace{.1in}
{\bf Proof.} Observe that if $f$ is finitely renormalizable, then it follows
immediately from Corollary 4.1. So, a new information we will get only
in infinitely renormalizable case.

Let $F_0\equiv I,\; l_1 $ be the first moment when orb(x) passes through $F_0$.
Now define inductively $F_k=[x(l_k), \tau(x(l_k))]$,
and $l_{k+1}$ as the first moment when orb($x$) passes
through the interval $F_k$. 

Denote by $M^-_k$ the component of cl\,($F_{k-1}\ssm F_k$) containing
$x(l_k)$, and by $M^+_k$ the component of  cl\,($F_k\ssm F_{k+1}$)
containing $x(l_k)$. Clearly, it is enough to show that for any $k$
there is an interval $J_k$ such that $M^+_k\subset J_k\subset F_k$ and
dens$(X|J_k)>1-\epsilon$.

Denote $o_k$=ord$(n,x,F_k)$.

If $o_k<2d$ then we have such an interval by Lemma 4.2. Otherwise 
by Proposition 2.2   there exists an interval $H_{k-1}\ni x$
monotonously mapped by $f^{l_k}$ onto $M^-_k\cup M^+_k$. Clearly,
dp$(n,x,M^-_k\cup M^+_k)\leq 2$, and we can apply the Two Interval Distortion
Lemma to $f^n|H_n$.

If $\lambda(M^+_k)\leq \lambda(M^-_k)$ then we conclude that $X$ is thick in
$M^+_k$ as required.
Otherwise $X$ is thick in $M^-_k$. Now let us pass from $x(l_{k})$
to $x(l_{k+1})$. It follows from Proposition 2.2 that 
ord($l_{k+1}-l_k,x(l_k), F_k)\leq d$. Hence, 
the existence of the required interval
$J_k$ follows again from Lemma 4.2.     $~~\Box$

\section{ Proof of the main results.}

Let us start with the following theorem proved by Guckenheimer [G] in the
case of negative Schwarzian derivative, and by Ma\~{n}\'{e} [Ma] in the general
smooth case (see also [vS]).

\vspace{.1in}
{\bf Theorem B.} 
{\em  For almost every $x\in \Lambda(f)$ the limit set $\omega (x)$
contains a critical point.}

\vspace{.1in}
{\bf Problem.} Is it true that
 $\omega(x)$ contains an {\em extremum} for almost every $x\in\Lambda(f)$?

\vspace{.1in}
{\bf Proof of Theorem 1.} Denote by $X_c$ the set $\{x: \omega (x)\ni c\}$.
By Theorem~A, 
\[\bigcup_{c\in C(f)} X_c =M \,({\rm mod}\, 0)\]
Let $X$ be any completely invariant subset of $X_c$ of positive measure.
In the finitely renormalizable case we have by Corollary 4.1 \newline
(i) If $c$ is an extremum then dens$X|c=1$ ;\newline
(ii) If $c$ is a reflection point then dens$X|c\geq 1/2$.
 
Hence, in the first case $f|X_c$ is ergodic, and in the second case
$X_c$ contains at most two ergodic components.

In the infinitely renormalizable case we can select $c$ as an extremum
and apply Lemma 4.4.  $~~\Box$

\vspace{.1in}
{\bf Proof of Theorem 2. } Let us associate to any ergodic component
$E\subset\Lambda (f) $ the following attractor $A$ (cf [M]).
 Consider the family of neighborhoods
$U$ such that the orbits of almost all points $x\in E$ pass through $U$
only finitely many times. Let $\{U_i\}_{i=1}^\infty$ be a countable basis of
this family. Set
\[A=M\ssm (\cup_{i=1}^\infty U_i).\]
Clearly, $A$ is a closed invariant set.
Let us check that   
\begin{equation}
 E\subset {\rm RL}(A)\; ({\rm mod}\, 0),
\end{equation}
Indeed, let $V$ be a neighbourhood of a point $a\in A$.
Let us consider the set

\vspace{.1in}
$E_V =\{x\in E : {\rm orb}(x)$ passes through $V$ infinitely many times\}.

\vspace{.1in}
\noindent
By definition, $\lambda (E_V)>0$. Since $E_V$ is completely invariant,
ergodicity yields $E_V = E\;({\rm mod}\, 0)$. Taking a countable basis of 
neighbourhoods $V_i$ of $A$, we obtain: $\omega(x)=A$ for almost all
$x\in E$ which is equivalent to (5.1).

Inclusion (5.1) implies that $A$ is a primitive attractor. 
So, we have constructed finitely many primitive attractors $A_k$ corresponding
to ergodic components $E_k$. Since
\[{\rm rl}(\cup A_k)\supset\cup E_k\supset\Lambda (f)\, \;({\rm mod\, 0}), \]
$\cup A_k$ is a global attractor for $f|\Lambda (f)$, 
and we have the decomposition (1.1).

Statements (i) and (ii) are also clear now: the former follows from (5.1),
the latter from Theorem B. Let us prove (iii).

Let $E_i$ and $E_j$ be two ergodic components in $\Lambda (f)$,
$A_i$ and $A_j$ be the corresponding primitive attractors. Assume
$A_i\cap A_j$ is infinite. Then we can pick up two close points
$a_1, a_2\in A_i\cap A_j$. Let us consider an interval $I=[a_0, a_3]$
containing $[a_1, a_2]$ and such that
\[|a_0-a_1|=|a_1-a_2|=|a_2-a_3|.\]
Applying Corollary 4.1, it is easy to see that each set $E_i,\; E_j$ is thick
at least in two of three intervals $ [a_0,a_1], [a_1, a_2], [a_2, a_3]$.
Hence, both $E_i, E_j$ are thick in one of these intervals, which is impossible.

Now let us pass to (iv). By (5.1), 
\[E_k\subset{\rm RL}(A_k)\subset {\rm rl}(A_k)\cap\Lambda(f).\]
So, it is enough to prove that the set
$({\rm rl}(A_k)\cap\Lambda(f))\ssm E_k$ has zero measure.
But otherwise there is
another ergodic component $E_i\subset\Lambda (f)\cap{\rm rl}(A_k ). $ 
The attractor $A_i$ corresponding to $E_i$ is infinite and is contained in
$A_k$, contradicting to what has been proved above.   $~~\Box$

\vspace{.1in}
{\bf Proof of Theorem 3.} It follows from Corollary 4.1 in the same way as in
[BL5].      $~~\Box$

\vspace{.1in}
{\bf Proof of Theorem 4.} (i) Let $K$ be a closed invariant set in $A$.
We are going to prove that $K=A$. Assume it is not the case. Then there is
a closed interval $I$ centered at a point $a\in A$ and such that 
$I\cap K=\emptyset$, $\lambda(I)<\gamma$ where $\gamma$ is taken from
Lemma 4.3.  Let $a$ divide $I$ into semi-intervals $I^+$ and $I^-$.
 Since $A$ is a Cantor attractor, there is an invariant set
$X\subset{\rm RL}(A)$ of positive measure such that
\begin{equation}
{\rm dens}(X|I^+)< 1-\epsilon\qquad{\rm and}\qquad{\rm dens}(X|I^-)
< 1-\epsilon.
\end{equation}

Let $N$ be from Lemma 4.3.
Find a neighborhood $U$ of $K$ such that orb$_N(x)\cap I=\emptyset$ for any
$x\in U$.
Let us pick up a density point $x\in X$ such that $\omega(x)=A$. Then
$f^nx\in U$ for some $n$, hence $I$ omits $N$ subsequent points
of orb(x), and we can apply Lemma 4.3. But its conclusion
contradicts (5.2). 
 
(ii) follows from Corollary 4.1. in the same way as in [BL5].

(iii) is immediate from (i) and Theorem 2(ii).

(iv) follows from (i) and Theorem 2(iii).   $~~\Box$

\vspace{.1in}
{\bf Proof of Theorem 5.} If cases A1 and A3 does not hold then
$A$ has non-empty interior. Since $A$ is primitive, it should be
 transitive. It is easy to show that a closed invariant transitive set
must be a cycle of intervals, so A2 holds. $~~\Box$

\vspace{.1in}
{\bf Proof of Corollary 1.1.} Associate to an infinite primitive attractor $A$
a critical point $c(A)$ belonging to $A$. Moreover, if $A$ is a cycle
of intervals, let $c(A)\in A^\circ$. It follows from Theorem 4(iv) that
this correspondence is injective, and we are done. $~~\Box$

\vspace{.1in}
{\bf Proof of Corollary 1.2.} It easily follows from Theorem 4(iv)
and \linebreak Theorem 5. $~~\Box$

\end{document}